\documentclass[11pt,a4paper]{amsart}

\usepackage[T1]{fontenc}
\usepackage{lmodern}
\usepackage{amsmath,amssymb,amsthm,mathtools}
\usepackage{mathrsfs}
\usepackage{geometry}
\usepackage{hyperref}

\theoremstyle{plain}
\newtheorem{thm}{Theorem}[section]
\newtheorem{lem}[thm]{Lemma}

\newtheorem{prop}[thm]{Proposition}

\theoremstyle{definition}

\newtheorem{rem}[thm]{Remark}

\newcommand{\R}{\mathbb R}

\newcommand{\Rk}{\R^{1+2d}}
\newcommand{\dd}{\, \mathrm{d} }
\newcommand{\dv}{\, \mathrm{d} v}
\newcommand{\dx}{\, \mathrm{d} x}
\newcommand{\dt}{\, \mathrm{d} t}

\newcommand{\D}{D_t}
\newcommand{\cS}{\mathcal S}
\newcommand{\cV}{V_{\alpha,\beta}}
\renewcommand{\L}{\operatorname{L}}

\title[Transfer from uncertainty]{Kinetic transfer of regularity  in any dimension \\  from uncertainty inequalities   on the real line}
\author{Pascal Auscher}
\address{Universit\'e Paris-Saclay, CNRS, Laboratoire de Math\'{e}matiques d'Orsay, 91405 Orsay, France}
\email{pascal.auscher@universite-paris-saclay.fr}
\date{November 22, 2026}

\thanks{We thank Lukas Niebel and Moritz Egert for helpful remarks and references. We used GPT 6-Astra only to find appropriate known literature on the one dimensional uncertainty inequality used here. The method to obtain transfer of regularity was developed solely by the author. A CC BY-NC 4.0 \url{https://creativecommons.org/licenses/by/4.0/} public copyright license has been applied by the authors to the present document and will be applied to all subsequent versions up to the Author Accepted Manuscript arising from this submission.
}
\subjclass[2010]{Primary: 42B10, 42B37, 35Q80 Secondary:   35Q20, 35Q84 }

\keywords{uncertainty inequalities, kinetic  transfer of regularity}

\begin{document}

\maketitle

\begin{abstract} The purpose of this short note is to show that one dimensional uncertainty inequalities easily induce transfer of regularity inequalities in kinetic equations in any dimension.
\end{abstract}

\section{Introduction}
 For kinetic equations with time $t\in \R$, position $x\in \R^d$ and velocity $v\in \R^d$, transfer of regularity means that  regularity along characteristics of a transport field $T=\partial_{t}+v\cdot\nabla_{x}$ and a velocity field $\nabla_{v}$ implies regularity along the position field $\nabla_{x}$. There could be also some transfer to regularity in time with possible loss in the velocity growth.  This was formulated in the celebrated hypoellipticity theory of H\"ormander \cite{MR0222474}, using  that the commutator between the two first fields yields the third one. Such phenomena admit very precise forms according to scaling rules. This is especially helpful for the regularity theory of the Landau and Boltzmann equations. 
 
 The connection between uncertainty inequalities and transfer of regularity in kinetic equations appears explicitly in an article of Alexandre et al \cite{AMUXY}. The point of that article  was to show that the kinetic transfer of regularity, after taking Fourier transform with respect to time $t$ and position $x$ and freezing the dual variables $\omega$ and $\varphi$, reduces to a Fefferman-Phong type inequality for certain operators with potential $q^2$, where $q=|\omega+v\cdot\varphi|$.  The leading term in  \cite[Theorem 1.1]{AMUXY} is   the fractional Laplacian $(-\Delta_{v})^\beta$, $0<\beta<1$.  The required Fefferman-Phong inequality when $\beta=1$  can be found in \cite[Theorem~6]{MM}. A special instance of the obtained transfer, when $0< \beta\le 1$, is
 \begin{equation}
\label{eq:AMUXY}
\|D_{x}^{\beta/(\beta+1)}f\|_{2}\lesssim  \|D_{v}^{\beta}f\|_{2}+ \|Tf\|_{2}
\end{equation} 
for  $f\in \L^2(\Rk)$ where $D_{x}=(-\Delta_{x})^{1/2}$, $D_{v}=(-\Delta_{v})^{1/2}$. It is homogeneity that provides the exact value of the transferred $x$-regularity exponent. We remark that this inequality, for all $\beta>0$ follows from Bouchut's inequality \cite[Proposition 1.1]{MR1949176}. The method there uses the same Fourier transform but does not rely on uncertainty estimates.

In a recent work  by Nystr\"om \cite{KN},  a similar reduction to a Fefferman-Phong inequality is used to prove
\begin{equation}
\label{eq:KN}
\|D_{x}^{1/3}f\|_{2}\lesssim  \|\nabla_{v}f\|_{2}+ \||T|^{1/2}f\|_{2}
\end{equation}
for  $f\in \L^2(\Rk)$ for which the right hand side is finite, where  $|T|$ is the positive self-adjoint square root of $-T^2$.  The  Schr\"odinger operator involved is $-\Delta_{v}+q$.    Observe that \eqref{eq:KN}  and  \eqref{eq:AMUXY} with $\beta=1$ share  the same homogeneity.

The questions we address is whether both inequalities above  belong to a larger family of transfer of regularity inequalities, whether they are all related to uncertainty estimates and whether they have a simple and common proof. Proofs of the uncertainty inequalities in the $v$-variable in \cite{AMUXY} use  pseudo-differential calculus combined with involved geometric arguments based on level sets of the potential.  The one in \cite{KN} is simpler because of the presence of the gradient but still has involved geometric reasoning.  Other tools using  dilation and  homogeneous group structure as the ones developed in \cite{KN}  could be quite delicate to implement for fractional exponents. 

The answer to those three questions is yes. The family of inequalities is the following statement.  

\begin{thm}
 For any $\alpha>0,\beta>0$, there a constant $C_{\alpha,\beta}>0$ such that for all $f$ in the space
\[
\mathcal{V}_{\alpha,\beta}=
\{
f\in \L^2(\Rk)\, ; \, D_{v}^\beta f\in \L^2(\Rk), |T|^\alpha f \in \L^2(\Rk)
\}
\]
 one has 
 \begin{equation}
\label{eq:KNalphabeta}
\|D_{x}^{\alpha\beta/( \beta+\alpha)}f\|_{2}\le C_{\alpha,\beta} ( \|D_{v}^\beta f\|_{2}+ \||T|^\alpha f\|_{2}).
\end{equation}
 \end{thm}
 
 When $0<\beta$ and $\alpha=1$, this is the same as  \eqref{eq:AMUXY} by unitary equivalence between $T$ and $|T|$.
 When $\beta=1$ and $\alpha=1/2$, this is the same as \eqref{eq:KN} by unitary equivalence between $\nabla_{v}$ and $D_{v}$. 
 
Let us come to a strategy of proof.  In \cite{AIN}, the inequality \eqref{eq:AMUXY} is proved  differently using the Galilean change of variables and partial Fourier transform in the position and velocity variables. The argument also uses  the fundamental solution of the Kolmogorov operator $(\partial_{t}+v\cdot\nabla_{x}) + (-\Delta_{v})^\beta$. In variables $\varphi$  dual to $x$ and $\xi$ dual to $v$,  the operator  becomes  $\partial_{t}+|\xi-t\varphi|^{2\beta}$, yielding the well-known expression to calculate the Fourier transform of the fundamental solution. So the argument there is designed to yield estimates involving $T$. This method was developed to obtain a  family of inequalities appropriate for weak solutions. Moreover, the specific use of $T$  gives access to other features  such as continuity in time in appropriate space for solutions of the equation and solvability of the Cauchy problem.  
 All of the  transformations used in \cite{AIN} being isometric on $\L^2(\Rk)$ (modulo a constant depending of the Fourier normalisation for Plancherel), it is natural to try this approach. 
 Rewriting \eqref{eq:KNalphabeta}
 in these new variables  will  allow us to see that this inequality easily   reduces  to a one dimensional uncertainty inequality for the operator $D_{t}^{2\alpha}+|t|^{2\beta}$ in $\L^2(\R)$. In other words, there is no geometric arguments as in \cite{AMUXY} and \cite{KN}. 
  
  This inequality is well-known. It appears in either multiplicative form 
  \begin{equation}
\|f\|_{\L^2(\R)}
\leq
M_{\alpha,\beta}
\big\||t|^\beta f\big\|_{\L^2(\R)}^{
\frac{\alpha}{\alpha+\beta}}
\big\|D_t^\alpha f\big\|_{\L^2(\R)}^{
\frac{\beta}{\alpha+\beta}},
\tag{M}
\end{equation}
or in equivalent additive form 
\begin{equation}
\|f\|_{\L^2(\R)}^2
\leq
A_{\alpha,\beta}
\left(
\big\||t|^\beta f\big\|_{\L^2(\R)}^2
+
\big\|D_t^\alpha f\big\|_{\L^2(\R)}^2
\right).
\tag{A}
\end{equation}
Both inequalities hold for every $f$ in the natural domain
\[
V_{\alpha,\beta}
=
\left\{
f\in \L^2(\R) \, ;\, 
|t|^\beta f\in \L^2(\R),\
D_t^\alpha f\in \L^2(\R)
\right\}.
\]
The multiplicative inequality can be obtained from the
entropic uncertainty principle of Hirschman \cite{Hirschman}. It is  a special case of Heisenberg-Pauli-Weyl inequalities by Cowling and Price \cite{CowlingPrice}. See \cite{FollandSitaram} for a survey. Such inequalities can be generalized to Riemannian manifolds and unimodular Lie groups equipped with  sublaplacians, see \cite{Martini} and references therein. Here again only the one dimensional case suffices.   Compactness
of the embedding of $V_{\alpha,\beta}$ into $L^2(\R)$ yields an easy proof of the additive inequality of Fefferman-Phong type.

In the sequel, we first  derive the transfer or regularity inequalities from (A). Next, we review two ways to obtain (A). One is  the entropic proof of (M)  and the other is a  proof of (A) by a self-contained compactness argument.  

\section{From uncertainty to transfer of regularity}

The Galilean change of variables defined by
\[\Gamma f(t,x,v)=f(t,x+tv,v)\]
is an isometry on $\L^2(\Rk)$. Our definition of (partial) Fourier transform in the variables $(\varphi,\xi)$ dual to $(x,v)$ is
\[
\widehat f(t,\varphi,\xi)=\iint_{\R^d\times \R^d} e^{-i(x\cdot \varphi+v\cdot\xi)}f(t,x,v) \dx \dv.
\]
For $f\in \L^2(\Rk)$, 
\[g(t,\varphi,\xi)=\widehat{\Gamma f}(t,\varphi,\xi)= \widehat f (t,\varphi, \xi-t\varphi),\]
so that $|T|^{\alpha}f$, $D_{v}^\beta f$ and $D_{x}^\gamma f$ respectively become $\D^{\alpha}g$,  $|\xi-t\varphi|^\beta g$ and  ${|\varphi|}^\gamma g$ for any $\gamma>0$, which are defined almost everywhere.  For the first term, we need to assume $|T|^{\alpha}f \in \L^2(\Rk)$.  Set $\gamma=  \frac{\alpha\beta}{ \beta+\alpha}$. 

Assume therefore $f\in \mathcal{V}_{\alpha,\beta}$.    By Plancherel's identity \eqref{eq:KNalphabeta} is equivalent to (with a comparable constant in a factor 2)  
\begin{equation}
\label{eq:KNalphabeta'}
  \int_{\Rk} \big(|\D^{\alpha}g|^2 + |\xi-t\varphi|^{2\beta}  |g|^2 \big)\dd(t,\varphi,\xi) \ge { \widetilde C_{\alpha,\beta}}^{-1}\int_{\Rk} {|\varphi|}^{2\gamma}| g|^2 \dd(t,\varphi,\xi).
\end{equation}
Invoking Fubini's Theorem, it suffices to show 
\begin{equation}
\label{eq:KNalphabeta't}
  \int_{\R} \big(|\D^{\alpha}u|^2 + |\xi-t\varphi|^{2\beta}  |u|^2 \big)\dd t \ge { \widetilde C_{\alpha,\beta}}^{-1}\int_{\R} {|\varphi|}^{2\gamma}| g|^2 \dd t
\end{equation}
for almost every $(\varphi,\xi)$ and $u(t)=g(t,\varphi,\xi)$. 

Fix one such  $(\varphi,\xi)$ with $\varphi\ne 0$, $u\in \L^2(\R)$ and the $t$-integral on the left hand side is finite.   
 Set $\delta=d(\xi, \R\varphi)$ the Euclidean distance of $\xi$ to the line $\R \varphi$.   Observe that 
\[
|\xi-t\varphi|^{2}=\delta^2+|t-t_{\varphi,\xi}|^2|\varphi|^2\ge |t-t_{\varphi,\xi}|^2|\varphi|^2 \] where $t_{\varphi,\xi}$ is the value of $t$ at which the second order polynomial $|\xi-t\varphi|^{2}$ realizes its minimum $\delta^2$. More precisely, $t_{\varphi,\xi}\varphi$ is the orthogonal projection of $\xi$ on $\R\varphi$ so that $\xi- t_{\varphi,\xi}\varphi$ is orthogonal to $(t-t_{\varphi,\xi})\varphi$ and the equality follows.
Hence,
\begin{align*}
 \int_{\R} \big(|\D^{\alpha}u|^2 + |\xi-t\varphi|^{2\beta}  |u|^2 \big)\dt
 \ge  \int_{\R} \big(|\D^{\alpha}u|^2  +  |t-t_{\varphi,\xi}|^{2\beta}|\varphi|^{2\beta}  |u|^2 \big)\dt.
\end{align*}
By making the change of variable $t\to \tau= t-t_{\varphi,\xi}$, we can replace $ t_{\varphi,\xi}$ by 0 and in particular $u(\tau+t_{\varphi,\xi})\in \cV$. Then applying (A) to $f_\lambda(\tau)=\lambda^{1/2}f(\lambda\tau)$ with $\lambda= |\varphi|^{-2\beta/\beta+\alpha}$, tracking scaling exponents and reversing the change of variable  yields that 
\begin{align*}
  \int_{\R} \big(|\D^{\alpha}u|^2  +  |t-t_{\varphi,\xi}|^{2\beta}|\varphi|^{2\beta}  |u|^2 \big)\dt
 \ge A_{\alpha,\beta}^{-1} |\varphi|^{2\gamma} \int_{\R}  |u|^2 \dt
\end{align*}
where $A_{\alpha,\beta}$ is the constant in (A).
 We have shown \eqref{eq:KNalphabeta't} with $\widetilde C_{\alpha,\beta}=A_{\alpha,\beta} $.
 Note that we have eliminated a term which is already controlled by the $D_{v}^\beta $ term alone.   

\section{The entropic proof of the multiplicative inequality}

We use  the following convention for the Fourier transform:
\[
\widehat f(\omega)
=
\int_{\R}e^{-it\omega}f(t)\, \dt,
\]
with inverse transform
\[
f(t)
=
\frac{1}{2\pi}
\int_{\R}e^{it\omega}\widehat f(\omega)\,\dd\omega.
\]
With this convention, Plancherel's identity reads in $\L^2(\R)$ norms
\[
\|f\|_{2}^2
=
\frac{1}{2\pi}
\|\widehat f\, \|_{2}^2.
\]

We first recall the one-dimensional entropic uncertainty
inequality of Hirschman.

\begin{thm}
Let $f\in\mathcal S(\R)$ satisfy $\|f\|_2=1$. Then
\begin{align}
-\int_{\R}|f(t)|^2\log |f(t)|^2\,\dt
-\int_{\R}
\frac{|\widehat f(\omega)|^2}{2\pi}
\log\left(
\frac{|\widehat f(\omega)|^2}{2\pi}
\right)\,\dd\omega
\geq \log(\pi e).
\tag{1.1}
\end{align}
\end{thm}

We combine it with the following  classical inequality for entropy by moment estimates, see, in particular,
\cite{LutwakYangZhang}.
The sharp constant and the equality cases are recalled
here for completeness.

\begin{lem}
Let $s>0$ and let $\rho$ be a probability density on
$\R$ such that
\[
m_s=\int_{\R}|x|^s\rho(x)\,\dx<\infty.
\]
Then
\begin{equation}
-\int_{\R}\rho(x)\log\rho(x)\,\dx
\leq
\log\left(K_s m_s^{1/s}\right),
\tag{1.2}
\end{equation}
where
\[
K_s
=
2e^{1/s}
\Gamma\left(1+\frac1s\right)s^{1/s}.
\]
\end{lem}

\begin{proof}
For $\lambda>0$, consider the probability density
\[
q_\lambda(x)
=
\frac{1}{Z_\lambda}
e^{-\lambda|x|^s},
\]
where
\[
Z_\lambda
=
\int_{\R}e^{-\lambda|x|^s}\,dx
=
2\Gamma\left(1+\frac1s\right)\lambda^{-1/s}.
\]
Integrating the inequality $\rho(x)\log\frac{\rho(x)}{q_\lambda(x)}\ge  \rho(x)-q_{\lambda}(x)$ yields 
\[
\int_{\R}\rho(x)\log\frac{\rho(x)}{q_\lambda(x)}\,dx
\geq 0.
\]
Therefore
\[
-\int_{\R}\rho(x)\log\rho(x)\,\dx
\leq
\log Z_\lambda+\lambda m_s.
\]
Choosing
\[
\lambda=\frac{1}{s m_s}
\]
yields
\[
-\int_{\R}\rho(x)\log\rho(x)\,\dx
\leq
\log\left(
2\Gamma\left(1+\frac1s\right)
(s m_s)^{1/s}
\right)+\frac1s.
\]
This is precisely (1.2). It is attained for that $q_{\lambda}$.
\end{proof}
Let us recast this towards (M).  For normalized functions, notice that $|f(t)|^2$ and $\frac{|\widehat f(\omega)|^2}{2\pi}$ are  probability densities. Apply the lemma  to the spatial density $|f(t)|^2$ with $s=2\beta$. We obtain
\begin{align}
-\int_{\R}|f(t)|^2\log|f(t)|^2\,\dt
&\leq
\log\left(
K_{2\beta}
\left(
\int_{\R}|t|^{2\beta}|f(t)|^2\,\dt
\right)^{1/(2\beta)}
\right)
\nonumber\\
&=
\log K_{2\beta}
+
\frac1\beta
\log\big\||t|^\beta f\big\|_2.
\tag{1.3}
\end{align}

For the frequency density $\frac{|\widehat f(\omega)|^2}{2\pi}$, take $s=2\alpha$. By Plancherel's
identity and the definition of $D_t^\alpha$,
\[
\int_{\R}|\omega|^{2\alpha}
\frac{|\widehat f(\omega)|^2}{2\pi}\,\dd\omega
=
\frac{1}{2\pi}
\int_{\R}|\omega|^{2\alpha}|\widehat f(\omega)|^2\,\dd\omega\\
=
\|D_t^\alpha f\|_2^2.
\]
Consequently,
\begin{align}
-\int_{\R}
\frac{|\widehat f(\omega)|^2}{2\pi}
\log\left(
\frac{|\widehat f(\omega)|^2}{2\pi}
\right)\,d\omega
\log K_{2\alpha}
+
\frac1\alpha\log\|D_t^\alpha f\|_2.
\tag{1.4}
\end{align}

Combining (1.1), (1.3), and (1.4), we obtain
\[
\frac1\alpha
\log\|D_t^\alpha f\|_2
+ \frac1\beta
\log\big\||t|^\beta f\big\|_2
\geq
\log(\pi e)
-\log K_{2\alpha}
-\log K_{2\beta}.
\]
Exponentiating gives
\[
\|D_t^\alpha f\|_2^{1/\alpha} \big\||t|^\beta f\big\|_2^{1/\beta}
\geq
\frac{\pi e}{K_{2\alpha}K_{2\beta}}.
\]
Raising both sides to the power
$\alpha\beta/(\alpha+\beta)$, we find
\[
\|D_t^\alpha f\|_2^{
\frac{\beta}{\alpha+\beta}}
\big\||t|^\beta f\big\|_2^{
\frac{\alpha}{\alpha+\beta}}
\geq
\left(
\frac{\pi e}{K_{2\alpha}K_{2\beta}}
\right)^{
\frac{\alpha\beta}{\alpha+\beta}}.
\]
This proves (M) for normalized functions. The general case
follows by homogeneity.

\section{The additive inequality by compactness}

We now give a proof of (A) based on compactness.
We equip $V_{\alpha,\beta}$ with the norm
\[
\|f\|_{V_{\alpha,\beta}}
=
\|f\|_2
+
\|D_t^\alpha f\|_2
+
\big\||t|^\beta f\big\|_2
.
\]
Remark that $\cS(\R)$ is easily shown to be densely contained in $V_{\alpha,\beta}$. 
\begin{prop}
For every $\alpha,\beta>0$, the embedding
\[
V_{\alpha,\beta}\hookrightarrow \L^2(\R)
\]
is compact.
\end{prop}

\begin{proof}  
Let $(f_n)$ be a bounded sequence in
$V_{\alpha,\beta}$. To establish that $(f_n)$ has a
subsequence converging strongly in $\L^2(\R)$,
we verify the two conditions (2.1) and (2.2) below of the Kolmogorov--M.~Riesz-Fr\'echet criterion (see, e.g., \cite{Brezis}, chapter 4).
First, for every $R>0$,
\[
\begin{aligned}
\int_{|t|>R}|f_n(t)|^2\,\dt
\leq
R^{-2\beta}
\int_{\R}|t|^{2\beta}|f_n(t)|^2\,\dt
=
R^{-2\beta}
\big\||t|^\beta f_n\big\|_2^2.
\end{aligned}
\]
Since $(f_n)$ is bounded in $V_{\alpha,\beta}$,
\[
\sup_n\int_{|t|>R}|f_n(t)|^2\,\dt
\longrightarrow0
\qquad\text{as }R\longrightarrow\infty.
\tag{2.1}
\]
Next, set
$
s=\min(\alpha,1)>0.
$
By Plancherel's identity,
\[
\begin{aligned}
\|f_n(\cdot+h)-f_n\|_2^2
&=
\frac1{2\pi}
\int_{\R}
|e^{ih\omega}-1|^2
|\widehat f_n(\omega)|^2\,\dd\omega.
\end{aligned}
\]
Using
\[
|e^{ih\omega}-1|^2
\leq C_s|h|^{2s}|\omega|^{2s} \le C_s|h|^{2s}(1+|\omega|^{2\alpha})
\]
we obtain
\[
\|f_n(\cdot+h)-f_n\|_2^2
\leq
C_s|h|^{2s}
\left(
\|f_n\|_2^2+\|D_t^\alpha f_n\|_2^2
\right).
\]
Thus
\[
\sup_n\|f_n(\cdot+h)-f_n\|_2
\longrightarrow0
\qquad\text{as }h\longrightarrow0.
\tag{2.2}
\]
\end{proof}

\begin{proof}[Proof of (A)]
We argue by contradiction. Suppose that the inequality
is false. Then there exists a sequence
$(f_n)\subset V_{\alpha,\beta}$ such that
\[
\|f_n\|_2=1 
\]
and
\[
\|D_t^\alpha f_n\|_2^2
+
\big\||t|^\beta f_n\big\|_2^2 
\longrightarrow0.
\tag{2.3}
\]
In particular, $(f_n)$ is bounded in
$V_{\alpha,\beta}$. By the compactness of the embedding,
after passing to a subsequence, there exists
$f\in \L^2(\R)$ such that
\[
f_n\longrightarrow f 
\qquad\text{strongly in }\L^2(\R).
\]
Hence
\[
\|f\|_2
=
\lim_{n\to\infty}\|f_n\|_2
=1.
\tag{2.4}
\]
Multiplication by $|t|^\beta$ being a closed operator
on $L^2(\R)$, we have using (2.3) that  $ |t|^\beta f=0$. 
 almost everywhere. 
Since $|t|^\beta>0$ almost everywhere, this implies
$f=0$ almost everywhere. This contradicts (2.4).
The contradiction proves the additive inequality.
\end{proof}

\section{Concluding remarks}

\begin{rem}
 The additive inequality implies the multiplicative one by a scaling argument, applying (A) to $f_r(t)=r^{1/2}f(rt)$ and optimizing. Conversely, use the arithmetic-geometric mean inequality. 
\end{rem}

\begin{rem}
 Our argument for \eqref{eq:KNalphabeta} avoids the use of density of,  say, $\cS(\Rk)$ in $\mathcal{V}_{\alpha,\beta}$. This is,  however,  true. Several ways to see that are possible but all have technical arguments.
\end{rem}

\begin{rem}
 Note that  compactness of $\mathcal{V}_{\alpha,\beta}\hookrightarrow \L^2(\Rk)$ is not clear since the $\L^2(\Rk)$ norm does not appear in \eqref{eq:KNalphabeta}. Identifying the completion of $\L^2(\Rk)$ for the homogeneous right hand side is not simple task. 
\end{rem}

\begin{rem}
 To make the note short, we restricted to this set of inequalities but to establish the validity of  Sobolev versions with negative or positive regularity and also  $L^p$ versions is of interest and probably more complicated. Some such inequalities are already in \cite{MR1949176} and \cite{AN}. The techniques developed in \cite{AN} or in \cite{N} could be helpful. 
\end{rem}

\end{document}